\documentclass[12pt]{article}    
    
\usepackage{latexsym,amsfonts,amsmath,epsfig,tabularx,amsthm,dsfont}    
    
\theoremstyle{definition}     
    
\newtheorem{thm}{Theorem}    
\newtheorem{rem}{Remark}

\newtheorem*{cor}{Main Theorem}

\newcommand\reals{\mathbb{R}}    
\newcommand\R{{\mathbb{R}}}    
\renewcommand\natural{\mathbb{N}}    
\newcommand\nat{\mathbb{N}}   
\newcommand\N{\mathbb{N}} 
\newcommand\sphere{\mathbb{S}} 
\newcommand\eps{\varepsilon}    
\newcommand\dist{{\rm dist}} 
\newcommand\Lip{{\rm Lip}}   
\newcommand\eu{{\rm e}}

\newcommand{\e}{\varepsilon}    
\newcommand{\wt}{\widetilde}

\newcommand{\points}{\mathcal{P}} 
\newcommand{\E}{\mathbb{E}}   
\newcommand{\ind}[1]{\mathds{1}_{#1}}    
\newcommand\dint{{\rm d}}  
 
\newcommand{\abs}[1]{\left\vert #1 \right\vert}  
\newcommand{\norm}[1]{\left\Vert #1 \right\Vert} 
 
\newcommand{\C}[1]{\mathcal{C}_{d}^{#1}}

\newlength{\fixboxwidth}    
\title{The Curse of Dimensionality for 
Numerical Integration of Smooth Functions} 
 
\author{Aicke Hinrichs,  
Erich Novak\footnote{This 
author was partially supported by the DFG-Priority Program 1324.}\\ 
Mathematisches Institut, Universit\"at Jena\\   
Ernst-Abbe-Platz 2, 07743 Jena, Germany\\   
email: a.hinrichs@uni-jena.de,   
erich.novak@uni-jena.de \\ 
\qquad   
\\ 
Mario Ullrich\footnote{This author was supported by DFG GRK 1523  
and ERC Advanced Grant PTRELSS.}\\  
Dipartimento di Matematica, Universit\`a Roma Tre \\ 
Largo San Leonardo Murialdo 1, 00146 Roma, Italy \\ 
email: ullrich.mario@gmail.com\\   
\qquad   
\\ 
Henryk Wo\'zniakowski\footnote{This author was partially    
supported by the National Science   
Foundation.   
}\\   
Department of Computer Science, Columbia University,\\   
New York, NY 10027, USA, and\\   
Institute of Applied Mathematics, University of Warsaw\\   
ul. Banacha 2, 02-097 Warszawa, Poland\\   
email:\ henryk@cs.columbia.edu}   
  
\begin{document}    
 
\maketitle    
 
\newpage 
  
\begin{abstract}    
We prove the curse of dimensionality for multivariate integration of 
$C^r$ functions:  
The number of needed function values to achieve an error $\eps$ is  
larger than $c_r (1+\gamma)^d$ for $\eps\le\eps_0$,  
where $c_r,\gamma>0$. 
The proofs are based on volume estimates for $r=1$ together  
with smoothing by convolution.  
This allows us to obtain smooth fooling functions for $r>1$.  
\end{abstract}     
 
\bigskip 
 
 \medskip 
 
 
\section{Introduction}    
 
We study multivariate integration 
for different classes $F_d$ of smooth functions  
$f \colon \R^d \to \R$.  
Our emphasis is on large values of $d\in\N$. 
We want to approximate 
\begin{equation}   \label{int} 
S_d(f) = \int_{D_d} f(x) \, \dint x \quad \mbox{for} \quad  
f\in F_d 
\end{equation} 
up to some error $\eps >0$, where  
$D_d \subset \R^d$ has Lebesgue measure 1. 
The results in this paper hold for arbitrary sets $D_d$,  
the standard example of course is $D_d = [0,1]^d$.   
 
We consider (deterministic)  
algorithms that use only function values, and   
classes $F_d$ of functions bounded in  
absolute value by 1 and containing all constant functions 
$f(x)\equiv c$ with $|c|\le 1$.  
An algorithm that uses no function value at all must be a constant,  
$A_0(f)\equiv b$, and its error is at least 
$$ 
\max_{f\in F_d}|S_d(f)|=1. 
$$ 
We call this the initial error of the problem, it does not depend on $d$.  
Hence multivariate integration is well scaled and that is why 
we consider $\eps<1$.   
 
Let  
$n(\eps,F_d)$  
denote the minimal number of  
function values needed for this task in the worst case  
setting\footnote{We add that $n(\eps,F_d)$ is the  
information complexity of multivariate integration over $F_d$ 
and is proportional to the (total) complexity as long as $F_d$ 
is convex and symmetric. The last two assumptions are needed to 
guarantee that a linear algorithm is optimal  and its  
implementation cost is linear in $n(\eps,F_d)$.}. 
By the \emph{curse of dimensionality} we mean that  
$n(\eps,F_d)$ is exponentially large in $d$.  
That is, there are positive numbers $c$, $\eps_0$ and $\gamma$ such that 
\begin{equation}\label{curse} 
n(\eps,F_d) \ge c \, (1+\gamma)^d \quad  
\mbox{for all} \quad \eps \le \eps_0  \quad  
\mbox{and infinitely many} \quad d\in \natural.  
\end{equation} 
For many natural classes  $F_d$ the bound in~\eqref{curse} 
will hold for all  
$d\in\natural$. 
This applies in particular to the classes considered in this paper. 
 
There are many classes $F_d$ for which the curse of dimensionality has  
been proved, 
see~\cite{NW08,NW10} for such examples. However, it has \emph{not} been 
known if the curse of dimensionality occurs for probably the most  
natural class which is the unit ball  
of $r$ times continuously differentiable functions,  
\[ 
\C{r}  
=\{f\in C^r( \R^d ) \ | \ \ \|D^\beta f\|\le 1 \quad 
\mbox{for all} \quad  |\beta|\le r\}, 
\] 
where $\beta=(\beta_1,\beta_2,\dots,\beta_d)$, with  
non-negative integers $\beta_j$, $|\beta|=\sum_{j=1}^d\beta_j$, 
and $D^\beta$ denotes the operator of  
$\beta_j$ times differentiation  
with respect to the $j$th variable for $j=1,2,\dots d$.  
By $\|\cdot\|$ we mean the sup norm, 
$\|D^\beta f\|=\sup_{x\in \R^d}  |(D^\beta f)(x)|$.   
 
For $r=0$, we obviously have $n(\eps, \C{0})=\infty$ for all 
$\eps<1$ and all $d\in\natural$.  
Therefore from now on we always assume that $r\ge1$.  
For $r=1$, the curse of dimensionality for $\C{1}$ follows from the results of  
Sukharev~\cite{Su79}. Whether the curse holds for $r\ge2$ has been  
an open problem for many years.  
 
The class $\C{r}$ for $D_d=[0,1]^d$ (and functions  
and norms restricted to $D_d$) 
was already studied in 1959 by  
Bakhvalov~\cite{B59}, see also~\cite{No88}. 
He proved that there are two positive  
numbers $a_{d,r}$ and $A_{d,r}$ such that 
\begin{equation}\label{bak1959} 
a_{d,r}\,\eps^{-d/r}\le n(\eps,\C{r})\le 
A_{d,r}\,\eps^{-d/r} \quad 
\text{for all}  \ d\in\nat 
\ \text{and} \ \eps\in(0,1). 
\end{equation} 
This means that for a fixed $d$ and for $\eps$ tending to zero,  
we know that $n(\eps,\C{r})$ is of order $\eps^{-d/r}$  
and the exponent of $\eps^{-1}$  
grows linearly in $d$.  
Unfortunately,  
Bakhvalov's result does not allow us to 
conclude whether the curse of dimensionality holds for the class $\C{r}$.  
In fact, if we reverse the roles of $d$ and $\eps$, and consider a fixed $\eps$ 
and~$d$ tending to infinity, the bound~\eqref{bak1959} 
on $n(\eps,\C{r})$ is useless. 
We prove the following result and hereby solve Open Problem 1 from~\cite{NW08}: 
 
\vspace{3mm} 
 
\begin{cor}\label{cor:Lip} 
The curse of dimensionality holds for the classes  
\,$\C{r}$ with the  
\emph{super-exponential} lower bound 
\[ 
n(\eps,\C{r}) \,\ge\, c_r\,(1-\eps) \, d^{\,d /(2r+3)} 
\quad 
\text{for all}  \ d\in\nat 
\ \text{and} \ \eps\in(0,1),  
\]  
where $c_r\in(0,1]$  
depends only on $r$. 
\end{cor}  
 
\vspace{3mm} 
We also prove that the curse of dimensionality holds 
for even smaller classes of functions~$F_d$ for which the norms of  
arbitrary directional derivatives are bounded proportionally to $1/\sqrt{d}$. 
 
We now discuss how we obtain lower bounds on $n(\eps,F_d)$ for 
numerical integration defined on  
convex and symmetric classes $F_d$.  
The standard proof technique is to find a fooling function $f\in F_d$ that  
vanishes at the points  
$\points = \{ x_1,x_2,\dots,x_n \} $ at which we sample functions 
from~$F_d$, and the integral of $f$ is as large as possible. 
All algorithms that use function values 
at $x_j$'s must give the same approximation of  
the integral for $f$ and $-f$. 
Thus, each such algorithm makes an error of at least  
$|S_d(f)-S_d(-f)|/2=|S_d(f)|$ for one of the functions. 
That is why the 
integral of $f$ is a lower bound on the worst case error 
of all algorithms using function values at $x_j$'s. 
If, for all choices of $x_1,x_2, \dots,x_n$,  
there are functions $f\in F_d$ vanishing at $x_j$'s  
with integrals larger than  
$\eps$ then $n(\eps,F_d)\ge n$.  
 
We start with the fooling function 
\[ 
f_0 (x) = \min\left\{1, \frac{1}{\delta\sqrt{d}}\, 
{\rm dist}(x,\points_\delta)\right\}  
\quad \mbox{for all} \quad x\in\R^d,  
\] 
where  
\[ 
\points_\delta = \bigcup_{i=1}^n B_\delta^d(x_i) 
\] 
and $B_\delta^d(x_i)$ is the ball with center $x_i$ and radius  
$\delta\sqrt{d}$.  
The function $f_0$ is Lipschitz. By a suitable smoothing  
via convolution we construct a fooling function  
$f_r \in \C{r}$ with $f_r|_\points = 0$.  
 
\section{Preliminaries} 
 
In this section, we precisely define our problem. Let $F_d$  be a class  
of  
continuous functions $f: \R^d \to\reals$ 
such that $S_d(f)$, see~\eqref{int}, exists for every $f\in F_d$.  
We approximate the integral $S_d(f)$, $f\in F_d$,  
by algorithms 
$$ 
A_{n,d}(f)=\phi_{n,d}\bigl(f(x_1),f(x_2),\dots,f(x_n)\bigr),   
$$ 
where $x_j\in \R^d $ can be chosen adaptively and $\phi_{n,d}:\reals^n\to 
\reals$ is an arbitrary mapping. Adaption means that the selection of $x_j$ 
may depend on the already computed values $f(x_1),f(x_2),\dots,f(x_{j-1})$. 
The (worst case) error of the algorithm $A_{n,d}$ is defined as 
$$ 
e(A_{n,d})=\sup_{f\in F_d}|S_d(f)-A_{n,d}(f)|. 
$$ 
The minimal number of function values to guarantee that the error is 
at most $\eps$ is defined as 
$$ 
n(\eps,F_d)= 
\min\{\,n\in\natural\ \big| \ \ \exists\ A_{n,d}\ \ \mbox{such that}\ \  
e(A_{n,d})\le\eps\}. 
$$ 
Hence we minimize $n$ over all choices of adaptive sample points $x_j$ and 
mappings $\phi_{n,d}$. It is well known that, as long as the class $F_d$ 
is convex and symmetric,  
we may restrict the minimization of $n$ by 
considering only nonadaptive choices of $x_j$ and  linear mappings $\phi_{n,d}$. 
Furthermore, 
\begin{equation}              
 \label{useful} 
  n(\eps,F_d)= 
  \min\Big\{\,n \in\natural\ | \ \  
  \inf_{\points\subset \reals^d, \#\points=n}\  
  \sup_{f\in F_d, f|_\points = 0 }  
  |S_d(f)|\le \eps\Big\},  
\end{equation} 
see~\cite[Prop.~1.2.6]{No88} 
or \cite[Theorem 5.5.1]{TWW88}.  
In this paper, we always consider convex and symmetric $F_d$  
so that we can use the last formula for $n(\eps,F_d)$.  
For more details see, e.g., Chapter 4 in~\cite{NW08}.

As already mentioned, 
our lower bounds 
are based on a volume estimate of a neighborhood of  
certain sets in $\R^d$, see also \cite{HNW11}.  
In the following, we denote by $A_\delta$  
the $(\delta\sqrt{d})$-neighborhood of $A\subset\R^d$,  
which is defined by 
\begin{equation}\label{eq02} 
A_\delta = \bigl\{ x \in \R^d \mid \  \dist(x, A) \le \delta\sqrt{d} \bigr\}, 
\end{equation} 
where $\dist(x,A)=\inf_{a\in A}\|x-a\|_2$  
denotes the Euclidean distance of $x$ from $A$.  
 
Furthermore, we denote by $B_\delta^d(x)$ 
the $d$-dimensional  
ball with center $x\in\R^d$ and radius $\delta\sqrt{d}$, i.e.,  
\[ 
B_\delta^d(x)\,=\, \bigl\{y\in\R^d \mid \ 
\norm{x-y}_2\le\delta\sqrt{d}\bigr\}. 
\] 
 
We will need some standard volume estimates for Euclidean balls. 
Recall that the volume of a Euclidean ball of radius 1 is given by 
$$  
V_d = \frac{\pi^{d/2}}{\Gamma(1+d/2)}. 
$$ 
{}From Stirling's formula for the $\Gamma$ function, we have 
$$  
\Gamma(x+1)= \sqrt{2 \pi x}\, x^x\, \eu^{-x+\frac{\theta_x}{12x}} \quad  
\mbox{for all} \quad \ x>0, 
$$ 
where $\theta_x\in(0,1)$,  
see \cite[p. 257]{AS72}. 
This leads to the estimate 
$$ \Gamma(x+1) > \sqrt{2 \pi x} \left( \frac{x}{\eu} \right)^x 
\quad \mbox{for all} \quad x>0.  
$$ 
Combining this estimate with the volume formula  
for the ball, we obtain for all $d\in\natural$, 
\begin{equation} \label{eq03}   
 \lambda_d\bigl(B_\delta^d(x)\bigr)  
 = \bigl(\delta\sqrt{d}\bigr)^d \, V_d 
 < \bigl(\delta\sqrt{d}\bigr)^d \, 
                \frac{\left(\frac{2 \pi \eu}{d} \right)^{d/2}}{\sqrt{\pi d}} 
 = \frac{\left( \delta \sqrt{2 \pi \eu}\right)^{d}}{\sqrt{\pi d}} 
 <  \left( \delta \sqrt{2 \pi \eu}\right)^{d}, 
\end{equation}  
where $\lambda_d$ is the Lebesgue measure. 
The volume formula for the Euclidean unit ball  
also shows the recurrence relation 
$$  
\frac{V_{d-1}}{V_d} = \frac{d}{d-1} \, \frac{V_{d-3}}{V_{d-2}} 
\quad \mbox{for all} \quad d\ge4. 
$$ 
This easily implies  
$$  
\frac{2}{\sqrt{d}} \, \frac{V_{d-1}}{V_d} < \frac{2}{\sqrt{d-2}} \,  
\frac{V_{d-3}}{V_{d-2}} 
\quad \mbox{for all} \quad d\ge4.  
$$ 
The last inequality can  
be used in an inductive argument leading to 
\begin{equation} \label{eq13}  
\frac{2}{\sqrt{d}} \,  
\frac{V_{d-1}}{V_d} \le 1 \quad \text{for all} \quad d\ge2. 
\end{equation}  
This will be needed later. 
 
\section{Convolution} \label{s:conv} 
 
In this section, we fix $k\in \natural$ and 
study the convolution  
$$f_k := f\ast g_1\ast\ldots\ast g_k 
$$  
of a function $f$ defined on $\R^d$ with 
(normalized) indicator functions $g_j$. 
We are interested in  
properties of $f_k$ in terms of the properties of  
the initial function $f$. 
Recall that the convolution of two functions $f$ and~$g$ on $\R^d$  
is defined by  
\[ 
(f\ast g)(x) = \int_{\R^d} f(x-t)\, g(t) \,\dint t 
\quad \mbox{for all}\quad x\in\R^d. 
\] 
Fix a number $\delta>0$ and a sequence  
$(\alpha_j)_{j=1}^k$ with $\alpha_j>0$  such that 
\[ 
\sum_{j=1}^k \alpha_j \le 1.  
\] 
For example, we may take  
$\alpha_j=1/k$  for $j=1,2\,\dots,k$. 
For $j=1,\dots,k$, we define the ball   
$$  
B_j \,=\, \Bigl\{x\in\R^d \,\big|\ \ 
\Vert x\Vert_2\,\le\,\alpha_j\, \delta \sqrt{d}\Bigr\}   
$$ 
and the function $g_j\colon \R^d\to\R$ by  
\begin{equation}\label{eq:gk} 
g_j(x) \,=\, \frac{\ind{B_j}(x)}{\lambda_d(B_j)}  
\,=\, \frac{1}{\lambda_d(B_j)}\,\begin{cases} 
1 & \  \text{ if }\, x\in B_j,\\ 
0 & \  \text{ otherwise. } 
\end{cases} 
\end{equation} 
Thus, the convolution of a function $f$ with  $g_j$  
can be written as 
\[ 
(f\ast g_j)(x) \,=\, \frac{1}{\lambda_d(B_j)}\,\int_{B_j}  
f(x+t)\,\dint t   \quad\mbox{for all}\quad x\in\R^d. 
\]  
We will frequently use the following probabilistic interpretation.  
Let $Y_j$ be a random variable that is uniformly distributed on $B_j$. 
Then the convolution of $f$ with $g_j$ can be written as the  
expected value 
\[ 
(f\ast g_j)(x) \,=\, \E\bigl[f(x+Y_j)\bigr]. 
\] 
 
The next theorem is the basis for the induction steps of the proofs  
of our main results.  
For $f\colon\R^d\to \R$, we use the 
Lipschitz constant 
\[ 
\Lip(f) = \sup_{x \not= y} \frac{|f(x)-f(y)|}{\Vert x-y\Vert_2}. 
\] 
Define 
\[ 
C^{r}  
= \bigl\{f\colon\R^d\to\R \mid\ \  
D^{\theta_\ell} \dots D^{\theta_1}f \text{ is continuous 
for all } \ell\le r \text{ and all\ }  
\theta_1,\dots,\theta_r\in \sphere^{d-1}\bigr\}, 
\] 
where $\sphere^{d-1}$ is the unit sphere in $\reals^d$ and  
$D^{\theta_1}f(x)=\lim_{h\to0}\frac1h 
\bigl(f(x+h\theta_1)-f(x)\bigr)$ is the derivative in the direction  
of $\theta_1$.  
 
\goodbreak 
 
\begin{thm}\label{thm:conv} 
For $k\in\N$ and $f\in C^r$, 
define  
$$ 
f_k=f\ast g_1\ast\ldots\ast g_k\quad 
\mbox{with\ \ $g_j$ from \eqref{eq:gk}}. 
$$ 
For $d\ge2$, let   
$\Omega \subset \R^d$ and let $\Omega_\delta$ be its  
neighborhood defined as in \eqref{eq02}. Then 
\begin{itemize} 
\item[$(i)$]   
if $f(x)=0$ for all  $x\in \Omega_\delta$ then  
$f_k(x)=0$ for all $x\in \Omega$,  
\item[$(ii)$] $\Lip(f_k)  \le  \Lip(f)$, 
\item[$(iii)$]  
if $\int_{D_d} f(x+t)\,\dint x \,\ge\, \eps$  
for all $t\in\R^d$ with $\Vert t\Vert_2\le\delta\sqrt{d}$ then  
$\int_{D_d} f_k(x) \dint x \ge \eps$,  
\item[$(iv)$] for all $\ell\le r$ and  
                all $\theta_1,\theta_2,\dots,\theta_\ell\in \sphere^{d-1}$,  
\[ 
\Lip\Bigl(D^{\theta_{\ell}}\,D^{\theta_{\ell-1}}  
\dots  D^{\theta_1}f_k\Bigr)  
\le \Lip\Bigl(D^{\theta_{\ell}}\,D^{\theta_{\ell-1}} 
 \dots  D^{\theta_1}f \Bigr),  
\] 
\item[$(v)$] $f_k\in C^{r+k}$,   
and for all $\ell\le r$, all $j=1,\dots,k$ and  
all $\theta_1,\theta_2,\dots,\theta_{\ell+j}\in \sphere^{d-1}$, 
\[ 
\Lip\Bigl(D^{\theta_{\ell+j}}\, D^{\theta_{\ell+j-1}}  
\dots  D^{\theta_{1}}f_k \Bigr)  
\le \biggl(\prod_{i=1}^j \frac{1}{\delta \alpha_{i}}\biggr) \, 
\Lip\Bigl(D^{\theta_{\ell}}\, D^{\theta_{\ell-1}}  
\dots  D^{\theta_{1}}f 
\Bigr) .  
\]  
\end{itemize} 
\end{thm} 
 
The parts $(i)$--$(iv)$  
of this theorem show  
that some properties  of the initial function~$f$ are 
preserved by convolutions. 
Part $(v)$ states that we gain one ``degree of smoothness'' 
with every convolution, losing only a multiplicative constant 
for its Lipschitz constant. 
 
\begin{proof} 
First note that 
we can write $f_k$ as  
\[ 
f_k(x) \,=\, \E\bigl[f(x+Y)\bigr], \quad \mbox{for all}\quad x\in\R^d, 
\] 
where $Y$ is a random variable with probability density function  
$g_1\ast\ldots\ast g_k$. By construction of $g_j$'s  which are the 
indicator functions of the balls whose sum of the radii is at most 
$\delta\sqrt{d}$, we have    
$$ 
\{t\in\R^d \mid \ g_1\ast\ldots\ast g_k(t)>0\}\subset 
\{t\in\R^d \mid \ \Vert t\Vert_2\le \delta\sqrt{d}\}, 
$$  
which implies that $x+Y\in\Omega_\delta$  
almost surely for every $x\in\Omega$.  
Thus, $f(x)=0$ for all $x\in \Omega_\delta$ implies  
that $f_k(x)=0$ for all $x\in \Omega$, which is property $(i)$. 
 
Property $(ii)$ is proven by  
\[\begin{split} 
\abs{f_k(x)-f_k(y)} \,&=\, \abs{\E\bigl[f(x+Y)-f(y+Y)\bigr]}  
\,\le\, \E\bigl[\abs{f(x+Y)-f(y+Y)}\bigr] \\ 
\,&\le\, \Lip(f)\; \E\bigl[\norm{(x+Y)-(y+Y)}_2\bigr] 
\,=\, \Lip(f) \norm{x-y}_2. 
\end{split}\] 
 
To prove $(iii)$, we use Fubini's theorem and we obtain 
\[ 
\int_{D_d} f_k(x)\,\dint x  
\,=\, \int_{D_d} \E\bigl[f(x+Y)\bigr]\,\dint x 
\,=\, \E\Bigl[\int_{D_d} f(x+Y)\,\dint x \Bigr]  
\,\ge\, \eps 
\] 
by assumption. 
 
For the proof of properties $(iv)$ and $(v)$, let  
$\theta=(\theta_1,\dots,\theta_\ell)\in(\sphere^{d-1})^\ell$. 
We write   
$D^\theta$ for $D^{\theta_{\ell}} \ldots D^{\theta_1}$.  
Clearly, $f\in C^{r}$ and $\ell\le r$ implies that $D^\theta f\in 
C^{r-\ell}\subseteq C$. Since $f_k$ is at least as smooth as $f$, 
both $D^\theta f$ and $D^\theta f_k$ are well defined.  
 
We need the well-known  
fact that $D^{\theta}(f\ast g)=(D^{\theta} f)\ast g$ if $f\in C^\ell$  
and $g$ has compact support.  
{}For 
$g=g_1\ast\ldots\ast g_k$,  
we have 
\[\begin{split} 
\abs{D^{\theta}f_k(x)- D^{\theta}f_k(y)}  
\,&=\, \abs{\bigl((D^{\theta}f)\ast g\bigr)(x)-  
\bigl((D^{\theta}f)\ast g\bigr)(y)}\\  
\,&=\, \abs{\int_{\R^d}\left[(D^{\theta}f(x+t) 
                        -D^{\theta}f(y+t)\right]\, g(t) \dint t} \\ 
\,&\le\, \Lip(D^\theta f)\, \norm{x-y}_2\, \int_{\R^d} g(t) \dint t\\ 
\,&=\, \Lip(D^\theta f)\, \norm{x-y}_2 
\end{split}\] 
for all $x,y\in \R^d$.   
The last equality follows since the $g_k$ is normalized.  
This proves $(iv)$. 
 
For $(v)$, we need to prove that $f_k\in C^{r+k}$ with $f_0=f\in 
  C^{r}$ by assumption, and then it is  
enough to show that for all $m\le r+k$  
and all $\theta=(\theta_{m},\dots,\theta_1)\in(\sphere^{d-1})^{m}$,  
\[ 
\Lip\Bigl(D^{\theta}f_{k} \Bigr)  
\le \frac{1}{\delta \alpha_{k}}\,\Lip\Bigl(D^{\bar\theta}f_{k-1} \Bigr),  
\] 
where $\bar\theta=(\theta_{m-1},\dots,\theta_1)\in(\sphere^{d-1})^{m-1}$.   
 
Assume inductively that $f_{k-1}\in C^{m-1}$, which holds for $k=1$. 
This implies  
$D^{\bar\theta}(f_{k-1}\ast g_k)=(D^{\bar\theta} f_{k-1})\ast g_k$, and 
\[\begin{split} 
D^\theta f_k(x) \,&=\,  
D^{\theta_{m}}\bigl((D^{\bar\theta}f_{k-1})\ast g_k\bigr)(x)\\ 
&=\, D^{\theta_{m}}\Bigl( \frac1{\lambda_d(B_k)} \int_{\R^d}  
                D^{\bar\theta}f_{k-1}(x+t)\, \ind{B_k}(t) \,\dint t\Bigr) \\ 
&=\, \frac1{\lambda_d(B_k)}\, D^{\theta_{m}}\Bigl(  
                \int_{\theta_{m}^\bot}  \int_{\R}  
                D^{\bar\theta}f_{k-1}(x+s+h\theta_{m})\,  
                \ind{B_k}(s+h\theta_m) \,\dint h \,\dint s \Bigr) \\ 
&=\, \frac1{\lambda_d(B_k)}\, \int_{\theta_{m}^\bot}\,  
                D^{\theta_{m}}\Bigl( \int_{\R}  
                D^{\bar\theta}f_{k-1}(x+s+h\theta_{m})\,  
                \ind{B_k}(s+h\theta_{m}) \,\dint h \Bigr)\,\dint s, \\ 
\end{split}\] 
where $\theta_{m}^\bot$ is the hyperplane orthogonal to $\theta_{m}$.  
For any function $f$ on~$\reals$ of the form 
\[  
f(x) = \int_{x-a}^{x+a} g(y)\, \dint y  
\] 
with some continuous function $g$ we have 
\[ 
f'(x) = g(x+a) - g(x-a). 
\] 
Therefore, we obtain 
\[\begin{split} 
D^\theta f_k(x) \,=\, \frac1{\lambda_d(B_k)}\,  
                \int_{B_k\cap\theta_{m}^\bot}\,  
&\biggl[ D^{\bar\theta}f_{k-1}\Bigl(x+s+h_{\rm max}(s)\,\theta_{m}\Bigr)\, \\ 
&\quad- D^{\bar\theta}f_{k-1}\Bigl(x+s-h_{\rm max}(s)\,\theta_{m}\Bigr) \biggr] 
\,\dint s 
\end{split}\] 
with 
\[ 
h_{\rm max}(s) \,=\, \max\{h\ge0 \mid \ \ s+h\theta_{m}\in B_k\}. 
\] 
{}For each $s\in{B_k\cap\theta_{m}^\bot}$, define the points  
$s_1=s+h_{\rm max}(s)\,\theta_{m}\in B_k$ and \\ 
$s_2=s-h_{\rm max}(s)\,\theta_{m}\in B_k$. Then 
\[\begin{split} 
\abs{D^{\theta}f_k(x)- D^{\theta}f_k(y)}  \,&\le\, \frac1{\lambda_d(B_k)}\,  
                \int_{B_k\cap\theta_{m}^\bot}\,\bigg[  
                \Bigl\vert D^{\bar\theta}f_{k-1}\bigl(x+s_1\bigr)\,  
                        - D^{\bar\theta}f_{k-1}\bigl(x+s_2\bigr) \\ 
&\qquad\qquad\qquad\qquad - D^{\bar\theta}f_{k-1}\bigl(y+s_1\bigr)\,  
                        + D^{\bar\theta}f_{k-1}\bigl(y+s_2\bigr)\Bigr\vert 
                \,\bigg]\dint s \\ 
&\le\, \frac1{\lambda_d(B_k)}\,  
                \int_{B_k\cap\theta_{m}^\bot}\,\bigg[ 
                \Bigl\vert D^{\bar\theta}f_{k-1}\bigl(x+s_1\bigr)\,  
                        - D^{\bar\theta}f_{k-1}\bigl(y+s_1\bigr) \Bigl\vert \\ 
&\qquad\qquad\qquad\qquad + \Bigl\vert D^{\bar\theta}f_{k-1}\bigl(x+s_2\bigr)\,  
        - D^{\bar\theta}f_{k-1}\bigl(y+s_2\bigr)\Bigr\vert \,\bigg]\dint s \\ 
&\le\, \frac{2\, \lambda_{d-1}(B_k\cap\theta_{m}^\bot)}{\lambda_d(B_k)}\,  
                \Lip(D^{\bar\theta}f_{k-1})\, \norm{x-y}_2. 
\end{split}\] 
In particular, this shows the implication  
\[ 
f_{k-1}\in C^{m-1} \,\Longrightarrow\, f_k\in C^m 
\]  
for all $k\in\N$.  
Taking $m=r+k$ we have $f_k\in C^{r+k}$, as claimed. 
 
For $m\le r+k$, it remains to bound  
$2\lambda_{d-1}(B_k\cap\theta_{m}^\bot)/\lambda_d(B_k)$.  
Recall that $B_k$ is a ball with radius $\delta\alpha_k \sqrt{d}$ and  
that $V_d$ is the volume of the Euclidean unit ball in $\reals^d$.  
We obtain from \eqref{eq13} that 
\[ 
\frac{2\, \lambda_{d-1}(B_k\cap\theta_{m}^\bot)}{\lambda_d(B_k)} 
\,=\, \frac{2 (\delta\alpha_k \sqrt{d})^{d-1}}{(\delta\alpha_k \sqrt{d})^{d}}\, 
                        \frac{V_{d-1}}{V_d}  
\,=\, \frac{2}{\delta\alpha_k \sqrt{d}}\, 
                        \frac{V_{d-1}}{V_d} 
\,\le\, \frac{1}{\delta\alpha_k}. 
\]  
This concludes the proof that  
\[ 
\Lip\Bigl(D^{\theta_{\ell+j}}\, D^{\theta_{\ell+j-1}}  
\dots  D^{\theta_{1}}f_k \Bigr)  
\le \biggl(\prod_{i=1}^j \frac{1}{\delta \alpha_{k+1-i}}\biggr) \, 
\Lip\Bigl(D^{\theta_{\ell}}\, D^{\theta_{\ell-1}}  
\dots  D^{\theta_{1}}f \Bigr),  
\]  
but since the order of convolution is arbitrary, we obtain in the same  
way 
\[ 
\Lip\Bigl(D^{\theta_{\ell+j}}\, D^{\theta_{\ell+j-1}}  
\dots  D^{\theta_{1}}f_k \Bigr)  
\le \biggl(\prod_{i\in J} \frac{1}{\delta \alpha_{i}}\biggr) \, 
\Lip\Bigl(D^{\theta_{\ell}}\, D^{\theta_{\ell-1}}  
\dots  D^{\theta_{1}}f \Bigr) 
\] 
for all $J\subset\{1,\dots,k\}$ with $\#J=j$.  
In particular, this implies $(v)$. 
\end{proof} 
 
\section{Main Results} 
 
Let $\points=\{x_1,\dots,x_n\}\subset\R^d$ 
be a collection of $n$ points.  
As pointed out in the introduction, we want to construct  
functions that vanish at $\points$ and have a large integral.  
For this, we choose  
\[ 
f_0(x) \,=\, \min\left\{1, \frac{1}{\delta\sqrt{d}}\, 
{\rm dist}(x,\points_\delta)\right\} \quad \mbox{for all}\quad x\in\R^d, 
\] 
where  
\[ 
\points_\delta \,=\, \bigcup_{i=1}^n B_\delta^d(x_i) 
\] 
and $B_\delta^d(x_i)$ is the ball with  
center $x_i$ and radius $\delta\sqrt{d}$.  
 
The function $\dist(\cdot,\points_\delta)$  
is Lipschitz with constant 1.  Hence, for $\delta\le1$, 
\begin{equation}\label{Lip} 
\Lip(f_0) = \frac{1}{\delta\sqrt{d}}. 
\end{equation} 
Additionally, $f_0(x)=0$ for all $x\in\points_\delta$ by definition.  
 
Using these facts we can apply Theorem~\ref{thm:conv}  
to prove the curse of dimensionality for the following class  
of functions that are defined on $\R^d$. 
For a fixed $r \in \nat$, we now take 
$\alpha_1=\dots=\alpha_r=\frac1r$ and define 
\[ 
F_{d,r,\delta} 
= \{f\colon \R^d \to \reals\ \big| \ \ f \in C^{r} \  
\mbox{satisfies~\eqref{cond1-Lip}--\eqref{cond3-Lip}}\}, 
\] 
where 
\begin{eqnarray}\label{cond1-Lip} 
 \Vert f \Vert &\le& 1,\\ 
\label{cond2-Lip} 
\Lip(f) &\le& \frac{1}{\delta\sqrt{d}},\\ 
\label{cond3-Lip} 
 \forall{k\le r}:\,\max_{\theta_1,\dots,\theta_k\in \sphere^{d-1}}\, 
        \Lip(D^{\theta_1} \dots  D^{\theta_k} f) &\le&  
        \frac{1}{\delta\sqrt{d}}\, \left( \frac{r}{\delta} \right)^k. 
\end{eqnarray} 
 
\begin{thm}\label{thm:Lip} 
For any $r\in\N$ and $\delta\in(0,1]$, 
\[ 
n(\eps,F_{d,r,\delta}) \,\ge\, (1-\eps) 
\begin{cases}1&\quad \mbox{for}\quad d=1,\\ 
\left( \delta\sqrt{18e\pi} \right)^{-d}& 
\quad \mbox{for}\quad d\ge2, 
\end{cases} 
\qquad 
\text{ for all } \eps\in(0,1). 
\]   
Hence the curse of dimensionality holds for  
the class $F_{d,r,\delta}$ for $\delta<1/\sqrt{18e\pi}$.  
\end{thm} 
 
This result shows that the growth rate of  
$n(\eps,F_{d,r,\delta})$ in $d$ can be arbitrarily large  
if we choose $\delta$ small enough. 
 
\begin{proof} 
Since the initial error for the classes $F_{d,r,\delta}$ is 1,  
we obtain $n(\eps,F_{d,r,\delta})\ge1$ for all $\eps\in(0,1)$. 
This proves the statement for $d=1$. 
 
For $d\ge2$,  
we use  
Theorem~\ref{thm:conv}  
with $k=r$,  
$\Omega=\points$ 
and  $f_{r}(x) \,=\, f_0\ast g_1\ast\ldots\ast g_{r}(x)$. 
Here, the $g_j$'s are as in Theorem~\ref{thm:conv}.  
Recall that we have chosen $\alpha_1=\ldots=\alpha_r=1/r$  
and $\alpha_j=0$ for $j>r$. 
The properties of the initial function $f_0$ and  
Theorem~\ref{thm:conv} immediately imply that 
$f_{r}$  satisfies~\eqref{cond1-Lip}--\eqref{cond3-Lip}.  
It remains to bound its integral. Note that 
$f_0(x)=1$ for all $x\notin\points_{2\delta}$.  
Clearly, $f_{r}(x)\ge0$ for all $x\in \R^d$.  
Since $f_{r}(x)$ depends only on the values  
$f_0(x+t)$ for $t\in\R^d$ with  
$\norm{t}_2\le\delta\sqrt{d}$, it follows that  
$f_{r}(x)=1$ for $x\notin\points_{3\delta}$.  
We thus obtain  
\[\begin{split} 
\int_{D_d}  f_{r}(x) \,\dint x  
\,&\ge\, \int_{D_d \setminus\points_{3\delta}} f_{r}(x) \,\dint x  
\,=\, 1- \lambda_d(\points_{3\delta}\cap D_d ) \\ 
\,&\ge\, 1-\lambda_d(\points_{3\delta})  
\,\ge\, 1- n\lambda_d(B_{3\delta}^d) \\ 
\,&>\, 1- \frac{n \Bigl(3\delta\sqrt{2e\pi}\Bigr)^d}{\sqrt{\pi d}}   \\ 
\,&>\, 1-n \Bigl(3\delta\sqrt{2e\pi}\Bigr)^d, 
\end{split}\] 
where the next to last inequality follows from the bound in  
\eqref{eq03}.  
Hence $\int_{D_d}f_{r}(x)\,\dint x \le \eps$ implies that 
$$ 
n\ge (1-\eps)\,(\delta\sqrt{18\e\pi})^{-d}. 
$$ 
Since this holds for arbitrary $\points$, the result follows. 
\end{proof} 
 
By Theorem~\ref{thm:Lip}, we  
know how the parameter $\delta$ comes into play.   
For $p>0$, let   
\[ 
\delta \,=\, \frac1{\sqrt{18e\pi}}\,d^{\,-p/(r+1)}. 
\] 
For this $\delta$, we obtain  
a somehow stronger form of the curse of dimensionality  
for the class 
\[ 
\wt{F}_{d,r,p} 
\,=\, \{f\colon \R^d \to \reals\ \big| \ \ f \in C^{r} \  
\mbox{satisfies~\eqref{cond1-Lip2}--\eqref{cond3-Lip2}}\}, 
\] 
where 
\begin{eqnarray}\label{cond1-Lip2} 
 \Vert f \Vert &\le& 1,\\ 
\label{cond2-Lip2} 
\Lip(f) &\le& d^{-\frac12+\frac{p}{r+1}}\;\sqrt{18e\pi},\\ 
\label{cond3-Lip2} 
 \forall{k\le r}:\,\max_{\theta_1,\dots,\theta_k\in \sphere^{d-1}}\, 
\Lip(D^{\theta_1} \dots  D^{\theta_k} f)  
&\le& d^{-\frac12+\frac{p (k+1)}{r+1}}\,r^k\,\bigl(\sqrt{18e\pi}\bigr)^{k+1}. 
\end{eqnarray} 
 
\begin{thm}\label{thm:Lip2} 
For any $r\in\N$ and $p > 0$,  
\[ 
n(\eps,\wt{F}_{d,r,p}) \,\ge\, (1-\eps)\, d^{\,p\,d /(r+1)} 
\quad 
\text{ for all }\;   d  \in \nat    
\;\text{ and }\;  \eps\in(0,1). 
\]   
Hence the curse of dimensionality holds for the class $\wt{F}_{d,r,p}$. 
\end{thm} 
 
Note that the classes $\wt{F}_{d,r,p}$ are contained in the classes 
\[ 
\C{r}=\{f\in C^r \ \big|\  \|D^\beta f\|\le 1 \quad  
\mbox{for all} \quad |\beta|\le r\}, 
\] 
if $p< 1/2$ and $d$ is large enough. 
This holds if 
\begin{equation}\label{eq:d0} 
d \,\ge\, \Bigl(r^{r}\,(18\eu\pi)^{(r+1)/2}\Bigr)^{1/(1/2-p)} .  
\end{equation} 
{}From this we easily obtain the main result  
already stated in the introduction.  
 
\vspace{2mm} 
 
\begin{cor} 
For any $r\in\N$, there exists a constant $c_r\in(0,1]$ 
such that 
\[ 
n(\eps, \C{r}) \ge  c_r\,(1-\eps)\, d^{\,d/(2r+3)} 
\quad 
\text{ for all }\; d\in\nat \;\text{ and }\;  \eps\in(0,1). 
\]   
Hence the curse of dimensionality holds for the class $\C{r}$. 
\end{cor} 
 
\vspace{2mm} 
 
\begin{proof} 
The case $d=1$ is trivial since the initial error for the classes  
$\mathcal{C}^r_d$ is again 1. 
 
For $d\ge2$, we  
know from Theorem~\ref{thm:Lip2} and the discussion thereafter  
that $n(\eps,\C{r})\ge (1-\eps)\, d^{\,p d /(r+1)}$  
for all $p<1/2$ if $d\ge d_0$, where $d_0=d_0(r,p)$ is the right  
hand side of \eqref{eq:d0}.  
This implies  
$$ 
n(\eps,\C{r})\ge \wt{c}_{r,p}\, (1-\eps)\,d^{\,p d/(r+1)} \quad 
\mbox{for all}\quad d\ge2. 
$$  
with 
$$ 
\wt{c}_{r,p}=d_0^{-pd_0/(r+1)}, 
$$ 
which depends only on $r$ and $p$. 
The choice $p^*=(r+1)/(2r+3)$ yields the result 
with $c_r=\wt{c}_{r,p^*}$. 
\end{proof} 
 
Note that  $c_r$ in the last theorem is super-exponentially small in $r$. 
 
\begin{rem} 
The reader might find it more natural to define classes  
of functions  
$F_{d,r}(D_d)$ that are defined only on $D_d \subset \R^d$. 
Not all such functions can be extended to smooth functions   
on $\R^d$, and even if they can be extended then the norm of the extended  
function could be much larger.  
Our lower bound results for functions defined on $\R^d$ 
can be also applied  
for functions defined on $D_d\subset\R^d$ 
and this makes them even stronger. 
\end{rem}  
 
\begin{rem}  
Note that the possibility of super-exponential lower bounds on the  
complexity depends on the definition of the  
Lipschitz constant.  
For the class  
$$ 
F_d = \left\{ f\colon [0,1]^d \to \R \mid \ \  
\sup_{x,y\in[0,1]^d}\frac{| f(x) - f(y)|}{\Vert x - y \Vert_\infty}\le 
1\right\},  
$$ 
Sukharev~\cite{Su79} proved that the product mid-point rule is optimal with  
error $e_n=\frac{d}{2d+2} n^{-1/d}$ for $n=m^d$.  
Hence, roughly, $n (\eps, F_d) \approx  
2^{-d} \eps^{-d} $ and the complexity is ``only'' 
exponential in $d$ for $\eps < 1/2$.  
\end{rem} 
 
\begin{rem}  
We mention two results for the very small class  
$$ 
F_d=C_d^\infty =\{f\in C^\infty ( [0,1]^d ) \mid\ \ \|D^\beta f\|\le 1 \quad  
\mbox{for all} \quad  \beta \in \N^d_0 \}. 
$$ 
O. Wojtaszczyk~\cite{Wo03} proved that  
$ 
\lim_{d \to \infty}  n(\eps, F_d) = \infty  
$  
for every $\eps <1$, hence the problem is not strongly polynomially tractable.  
It is still open whether the curse of dimensionality holds for  
this class $F_d$. 
The same class $F_d$ was studied  
for the approximation problem in \cite{NW09}. 
For this problem the curse of dimensionality is present 
even if we allow algorithms that use arbitrary linear functionals. 
\end{rem}

\medskip 
 
\noindent 
{\bf Acknowledgement.} 
We thank Jan Vyb\'iral and Shun Zhang for valuable remarks.

\end{document}